\documentclass[11pt, oneside]{article}   	

\usepackage[utf8]{inputenc} 
\usepackage{textcomp} 
\usepackage{graphicx}  
\usepackage{flafter}  
\usepackage{natbib} 

\usepackage{amsmath,amssymb}  
\usepackage{bm}  

\usepackage[pdftex,bookmarks,colorlinks,breaklinks]{hyperref}  
\hypersetup{linkcolor=red,citecolor=blue,filecolor=dullmagenta,urlcolor=darkblue} 

\usepackage{memhfixc}  
\usepackage{pdfsync}  

\usepackage{geometry}                		
\usepackage{graphicx}				
\usepackage{amssymb}

\newcommand{\bihomographic}[8]{\left(\begin{smallmatrix}#1&#2&#3&#4\\#5&#6&#7&#8\end{smallmatrix}\right)}

\newcommand{\abcd}{\left(
\begin{smallmatrix} 
a & b & c & d\\ 
e & f & g & h
\end{smallmatrix}
\right)}

\newcommand{\produce}{\mbox{prod}}

\title{Computing with Continued Logarithms}
\author{Michael J. Collins\\Daniel H. Wagner Associates\\mjcollins10@gmail.com}
\date{} 

\begin{document}
\maketitle
\begin{abstract}
Gosper  developed an algorithm for performing arithmetic on continued fractions (CFs), and introduced continued logarithms (CLs) as a variant of continued fractions better suited to representing extremely large (or small) numbers. CLs are also well-suited to efficient hardware implementation. Here we present the algorithm for arithmetic on CLs, then extend it to the novel contribution of this paper, an algorithm for computing trigonometric, exponential, and log functions on CLs. These methods can be extended to other transcendental functions.

As with the corresponding CF algorithms, computations are entirely in the domain of the CL representation, with no floating-point arithmetic; we read one CL input term at a time, producing the next CL term of the result as soon as it is determined. The CL algorithms are in fact simpler than their CF counterparts.

 We have implemented these algorithms in Haskell.
\end{abstract}

\section{Introduction and Motivation}
Let $x$ be  a real number greater than $1$. Following \cite{hakmem}, if $\lfloor \log_2(x)\rfloor = k$, i.e. if $2^k \leq x < 2^{k+1}$, we write
\begin{equation}\label{eq:CLdefinition}
x = 2^k\left(1 + \frac{1}{y}\right) 
\end{equation}
and recursively define the the \emph{continued logarithm} (CL) representation of $x$ to be the sequence of integers starting with $k$, followed by the CL representation of $y = \frac{1}{2^{-k}x-1}$; clearly $y > 1$. We can view the first term as the order of magnitude of $x$, followed by the inverse order of magnitude of the ``adjustment" needed to improve the initial approximation, and so on recursively. The CL representation of $2^k$ is just $[k]$. We may write a CL representation $x=[a_0,a_1,\cdots]$ as $x=[a_0,a_1,\cdots]_2$ when we need to explicitly distinguish continued \emph{logarithms}  from continued \emph{fractions}.

As an example, we find the CL representation of 19. Sixteen is the largest power  of $2$ less than $19$, so we start with
\[
19 = 2^4\frac{19}{16} = 2^4\left(1 + \frac{1}{\frac{16}{3}}\right) \ .
\]
Now $2^2 < \frac{16}{3} < 2^3$, so
\[
\frac{16}{3} = 2^2\left(1 + \frac{1}{3}\right)
\]
and
\[
3 = 2^1\left(1 + \frac{1}{2}\right) \ .
\]
We finish by representing 2 as $[1]$, so $19 = [4,2,1,1]$.

Similarly, the representation of $\pi$ starts with 1; then $\frac{1}{\pi/2 -1} = 1.751938\cdots$ and
$\frac{1}{0.751938\cdots} = 1.329896\cdots$, so the next two terms are $0,0$ and it continues
\begin{equation}\label{eq:CLpi}
\pi = [1,0,0,1,0,0,3,0,3,0,2,0,0,2,5,...]\ .
\end{equation}
In the same manner
\begin{equation}\label{eq:CLe}
e = [1,1,1,1,0,2,2,0,2,0,0,0,1,1,0,...]\ .
\end{equation}
In contrast to its continued fraction, the CL terms of $e$ have no simple pattern.
An especially interesting case is the golden ratio; from $\phi = \frac{1}{\phi-1}$ we obtain $\phi = [0,0,0,\cdots]$.

The basic properties of CLs are well-covered in \cite{Borwein02102017} (and see \cite{rockett1992continued} for continued fractions). As with CFs, rational numbers always have finite CLs, and periodic CLs represent quadratic irrationals. However, not all quadratic irrationals have periodic CL terms; it is an interesting open problem to characterize the integers $n$ for which $\sqrt{n}$ is CL-periodic. An analogue of Khinchin's theorem holds for CLs, and it is not surprising that CL terms tend to be small; the largest among the first ten thousand terms of $\pi$ and $e$ are 15 and 18.

\subsection{Dynamical System Representation}
There is another way of writing CLs which is trivially equivalent, but better suited to the computational algorithms of this paper.
We can define a dynamical system on $[1,\infty)$ by
\begin{equation}\label{eq:dynSys}
g(x) = \begin{cases}
             x/2 & \mbox{if } x \geq 2 \\
                 \frac{1}{x-1} & \mbox{if } 1 < x < 2 \\
                 \mbox{terminate} & \mbox{if } x = 1
              \end{cases}
\end{equation}
We write ``1" for every application of the first rule and ``0" for every application of the second, so, for instance, the CL of 19 becomes $[1,1,1,1,0,1,1,0,1,0,1]$; the more compact representation counts the number of ones immediately preceding each zero or termination.

\subsection{Representing $x < 1$}
There are at least two plausible ways to represent $0 < x < 1$. Since we know that $x^{-1}>1$, we could take an ``impatient" approach and simply set the initial term to $-1$, followed by the representation of $x^{-1}$. As a dynamical system this would be
\begin{equation}
g'(x) = \begin{cases}
             x/2 & \mbox{if } x \geq 2 \\
                 \frac{1}{x-1} & \mbox{if } 1 < x < 2 \\
                 \frac{1}{x} & \mbox{if } 0 < x < 1 \\
                 \mbox{terminate} & \mbox{if } x = 1
              \end{cases}
\end{equation}
writing $-1$ for application of the third rule. As with ordinary continued fractions, $x$ differs from $x^{-1}$ only by the addition or removal of an initial term.

On the other hand, it seems more in the spirit of CLs to multiply $x$ by $2$ as many times as we need to get $2^kx \geq 1$, so the more ``natural" representation would be
\begin{equation}
\hat{g}(x) = \begin{cases}
                 x/2 & \mbox{if } x \geq 2 \\
                 \frac{1}{x-1} & \mbox{if } 1 < x < 2 \\
                 2x & \mbox{if } 0 < x < 1 \\
                 \mbox{terminate} & \mbox{if } x = 1
              \end{cases}
\end{equation}
again writing $-1$ for application of the third rule. Note that in this case (\ref{eq:CLdefinition}) continues to hold for negative $k$. Now the representations of $x$ and $x^{-1}$ are not obviously related; for instance $\phi^{-1}$ becomes $[-1,2,4,4,4,4,\cdots]$, and $1/19$ is $[-5,0,0,1,3,1]$.

Our algorithms will use the first ``impatient" method for the sake of simplicity (following \cite{4402399,5467052}); the modifications needed for the second ``natural" method would be straightforward.

Negative numbers may be represented as $x=[-2,x_1,\cdots]$ where $-x=[x_1,x_2\cdots] > 0$.
The empty list implicitly represents infinity (since $2^k = 2^k(1 + \frac{1}{\infty})$),  so zero may be represented impatiently as $[-1]$, or naturally as $[-1,-1,-1,\cdots]$.

\section{Arithmetic on Continued Logarithms}\label{sec:arithAlg}
We now consider how to do arithmetic on CLs. The algorithm is implicit in \cite{hakmem}, and stated very concisely in \cite{5467052} (which describes a hardware implementation), but we hope our exposition here will be of some value; it will also allow us to introduce some notation we will need later.

How would we, for example, derive the continued logarithm $[z_0,z_1,\cdots] = \pi/e = 1.155727\cdots$ solely from (\ref{eq:CLpi}) and (\ref{eq:CLe})? As already noted, it will be convenient to use the (admittedly hard-to-read) ``dynamical system" representation which gives
\[
\pi = [1,0,0,0,1,0,0,0,1,1,1,\cdots] = [p_0, p_1,p_2,\cdots]
\]
and
\[
e = [1,0,1,0,1,0,1,0,0,1,1,0,\cdots] = [e_0,e_1,e_2,\cdots]
\]
with all $p_i, e_i \in \{0,1\}$. Since $\pi = 2(1+\frac{1}{[p_2,p_3\cdots]})$ and $e=2(1+\frac{1}{[e_2,e_3\cdots]})$, we can begin with
\[
\frac{\pi}{e} = \frac{1+\frac{1}{[p_2,p_3\cdots]}}{1+\frac{1}{[e_2,e_3\cdots]}}
     = \frac{[p_2,\cdots][e_2,\cdots] + [e_2,\cdots]}{[p_2,\cdots][e_2,\cdots] + [p_2,\cdots]}
\]
then substituting $[p_2,\cdots] = 1 + \frac{1}{[p_3,\cdots]}$ and $[e_2,\cdots] = 2[e_3,\cdots]$ takes us to
\[
\frac{\pi}{e} = \frac{4[p_2,\cdots][e_2,\cdots] + 2[e_2,\cdots]}{2[p_2,\cdots][e_2,\cdots] + [p_2,\cdots]+ 2[e_2,\cdots] + 1}
\]
We have now read enough terms to imply $e < \pi < 2e$, and it is easy to observe that the above expression must lie between 1 and 2 whenever $[p_2,\cdots]$ and $[e_2,\cdots]$ are both greater than 1; so the first term of $\pi/e$ is $z_0=0$, and
\[
[z_1,z_2,\cdots] = \left(\frac{4[p_2,\cdots][e_2,\cdots] + 2[e_2,\cdots]}{2[p_2,\cdots][e_2,\cdots] + [p_2,\cdots]+ 2[e_2,\cdots] + 1}  -1\right)^{-1} 
\]
i.e.
\[
[z_1,z_2,\cdots] =  \frac{2[p_2,\cdots][e_2,\cdots] + [p_2,\cdots]+ 2[e_2,\cdots] + 1}{2[p_2,\cdots][e_2,\cdots]  -[p_2,\cdots] - 1}
\].

This is not sufficient to determine $z_1$; the bihomographic expression ranges from one to infinity as $[p_2,\cdots]$ and $[e_2,\cdots]$ independently range from one to infinity.
So we must make further substitutions, beginning with $[p_2,\cdots] = [0,p_3,\cdots] = (1 + \frac{1}{[p_3,\cdots]})$ and so on.

The key idea of the arithmetic algorithm is that, as we continue to make such substitutions, the still-undetermined output terms $[z_j, z_{j+1},\cdots]$ are always given by an expression of the form
\begin{equation}\label{eq:bihomographic}
M(x,y) = \frac{axy + bx + cy + d}{exy + fx + gy + h}
\end{equation}
where $x = [x_k, x_{k+1},\cdots]$ and $y=[y_k, y_{k+1},\cdots]$ are the still-unread input terms, and all coefficients $a,b \ldots h$ are integers. Such an expression is called \emph{bihomographic}, and can be
identified with the matrix $M=\abcd$. It is straightforward\footnote{conceptually it is straightforward, but determining these bounds is most of the computational work of arithmetic; Brabec \cite{5467052} considers hardware optimizations.} 
to determine upper and lower bounds on $M(x,y)$ subject to the conditions $x \geq 1, y \geq 1$;
when these conditions imply $1 < M(x,y) < 2$, we know that $z_j=0$, and the expression for $[z_{j+1}, z_{j+2}\cdots]$ is
\[
\left(\frac{axy + bx + cy + d}{exy + fx + gy + h} - 1\right)^{-1} = \frac{exy + fx + gy + h}{(a-e)xy + (b-f)x + (c-g)y + (d-h)}
\]
so we define the corresponding ``produce output" operation on matrices
\[
\mbox{prod}_0\abcd = \bihomographic{e}{f}{g}{h}{(a-e)}{(b-f)}{(c-g)}{(d-h)}
\]
Similarly if $M(x,y) \geq 2$, we know that $z_j=1$, and the expression for $[z_{j+1}, z_{j+2}\cdots]$ is
\[
\frac{axy + bx + cy + d}{2exy + 2fx + 2gy + 2h}
\]
so we define 
\[
\mbox{prod}_1\abcd = \bihomographic{a}{b}{c}{d}{2e}{2f}{2g}{2h} \ .
\]

 If $z_j$ is not yet determined, we must narrow the range by consuming the next term of either $x$ or $y$.
 If $x = [0,x_{k+1},x_{k+2}\cdots]$, we can make the substitution $x \leftarrow 1 + 1/x$ to get
\[
[z_j, z_{j+1},\cdots] = \frac{(a+c)xy + (b+d)x + ay + b}{(e+g)xy + (f+h)x + ey + f}
\]
with $x = [x_{k+1},x_{k+2}\cdots]$. In terms of matrices we can define a ``consume input" operation,
\[
\mbox{con}_{x,0}\abcd = \bihomographic{a+c}{b+d}{a}{b}{e+g}{f+h}{e}{f} \ .
\]

Similarly $x = [1,x_{k+1},x_{k+2}\cdots]$ implies the substitution $x \leftarrow2x$ and 
\[
[z_j, z_{j+1},\cdots]  = \frac{2axy + 2bx + cy + d}{2exy + 2fx + gy + h}
\]
with 
\[
\mbox{con}_{x,1}\abcd = \bihomographic{2a}{2b}{c}{d}{2e}{2f}{g}{h} \ .
\]
If we reach the end of a finite input, we substitute $x \leftarrow 2$ or $x \leftarrow 1$ as $x=[1]$ or $[0]$. The substitutions for $y$ are analogous.

In particular, addition ($z = x+y$) and multiplication  ($z = xy$) are just special cases of computing $M(x,y)$, with
$M = \left(
\begin{smallmatrix} 
0 & 1 & 1 & 0\\ 
0 & 0 & 0 & 1
\end{smallmatrix}
\right)$ and 
$M = \left(
\begin{smallmatrix} 
1 & 0 & 0 & 0\\ 
0 & 0 & 0 & 1
\end{smallmatrix}
\right)$ respectively. Subtraction and division are $M=\bihomographic{0}{1}{-1}{0}{0}{0}{0}{1}$ and $M=\bihomographic{0}{1}{0}{0}{0}{0}{1}{0}$.

This is the basic notion of Gosper's arithmetic algorithm, but a full implementation must also consider the issue of termination.
The problem is that repeatedly consuming input is not guaranteed \emph{ever} to determine the next $z_j$.
As we read more terms of $x$ and $y$, we might just obtain increasingly tight bounds of the form $2-\varepsilon < M(x,y) < 2 + \varepsilon$, which will never tell us if $z_j$ is 1 or 0.
This is what must happen if $x$ and $y$ are infinite while $z$ is finite and $[z_j z_{j+1},\cdots] = [z_j] = [1]$; finite prefixes of $x$ and $y$ do not determine the result.
There are various ways to approach this problem; we will  put the issue aside for the moment and return to it in section \ref{sec:termination}. Our Haskell \cite{hutton2007programming} implementation\footnote{https://github.com/mjcollins10/ContinuedLogarithms} currently does not handle this case.

\section{Transcendental Functions on Continued Logarithms}\label{sec:transcendental}
We now extend the ideas of the arithmetic algorithm to compute $\sin, \arcsin, \exp$, and $\log$ over CLs. The approach is very similar to that used for continued fractions in \cite{collins2025transcendentalfunctionscontinuedfractions}.

\subsection{Exponentials} We first consider how to obtain the CL terms of 
\[
e^x = \sum_{k=0}^{\infty} \frac{x^k}{k!}
\]
where $x$ itself is given as a continued logarithm. We will assume $0 \leq x \leq 1$; otherwise we can use the identity $e^x = (e^{x/2})^2$. Let
\begin{equation}\label{eq:Mxn}
M^x_n(y) = 1 + \frac{xy}{n} \ .
\end{equation}
This is a bihomographic expression with matrix $M_n=\bihomographic{1}{0}{0}{n}{0}{0}{0}{n}$.
The $n^{th}$ degree Taylor approximation
$1 + x + x^2/2 +\cdots + x^n/{n!}$ can be written using Horner's rule as
\[
 1 + x \left(1 + \frac{x}{2}\left( 1 + \frac{x}{3}\left(1 + \cdots \left(1 + \frac{x}{n}\right)\right)\right) \cdots \right) = M^x_1( M^x_2( \cdots M^x_n(1) \cdots ) )
 \]
 so we define
 \begin{equation}\label{eq:yn}
 y_n(x) = 1 + \frac{x}{n} + \frac{x^2}{n(n+1)} + \cdots = \lim_{m\to\infty} M^x_n( M^x_{n+1}( \cdots M^x_m(1) \cdots ) ) \ .
\end{equation}
Now $y_n(x) = M^x_n(y_{n+1}(x))$, and $y_1(x) = e^x$.
Since $x\leq1$ we have
 \[
 y_n(x) \leq 1 + \frac{e^x - 1}{n} < 1 + \frac{2}{n}\ .
 \]
 For $n \geq 2$ this gives $y_n(x)<2$, so the first CL term of $y_n(x)$ is zero. Furthermore 
 \[
 \frac{1}{y_n(x) - 1} > \frac{n}{2}
  \]
 so the initial zero term is followed by $\lfloor \log_2(n)-1\rfloor$ ones. We get the subsequent terms of $y_n(x)$ by applying the arithmetic algorithm to
 \[
 \hat{M_n}(x, y_{n+1}(x))
  \]
 where $\hat{M_n}$ is $\mbox{prod}^{\lfloor \log_2(n)-1\rfloor}_1(\mbox{prod}_0 (M_n))$, with $\mbox{prod}_0, \mbox{prod}_1$ as defined in section \ref{sec:arithAlg}.
 The terms of $y_{n+1}(x)$ are obtained recursively in the same manner.

 \subsubsection{Termination}
 We have seen that we can apply the arithmetic algorithm to an infinitely nested expression like \ref{eq:yn},
 because we always know the first few terms of $y_n(x)$ without recursing into $y_{n+1}(x)$.
 It remains to show that we will in fact obtain every term of $y_n(x)$ eventually. 
 In other words we must show that, given any $\varepsilon > 0$,
 we eventually confine $y_n(x)$ to an interval of radius $\varepsilon$.
 From \ref{eq:Mxn} we see that we will have
\[
 K - \varepsilon < y_n(x) < K + \varepsilon
\]
for some $K$ if, for some $K'$,
 \[
 K' - \frac{n\varepsilon}{x} < y_{n+1}(x) < K' + \frac{n\varepsilon}{x}
\]
and thus if
 \[
 K'' - \frac{n(n+1)\cdots(n+r-1)\varepsilon}{x^r} < y_{n+r}(x) < K'' + \frac{n(n+1)\cdots(n+r-1)\varepsilon}{x^r}
 \]
for some $r$. Since $\frac{n(n+1)\cdots(n+r-1)}{x^r}$ approaches infinity,
we will eventually get the required bound from the known initial terms of $y_{n+r}(x)$,
and will not need to recurse into $y_{n+r+1}(x)$
(assuming we have also ingested enough terms of $x$ to get a sufficiently tight bound on $x$).

Similar observations apply to computation of the other transcendental functions we now consider.
 
 \subsection{Inverse Trigonometric Functions}
We can apply a similar approach to
\[
\frac{\arcsin(x)}{x} = \sum_{n=0}^{\infty} \frac{(2n)!}{4^n(n!)^2(2n+1)}x^{2n} = 1 + \frac{x^2}{6} + \frac{3x^4}{40} + \cdots
\]
with $| x | \leq 1$. Let $t_n = \frac{(2n)!}{4^n(n!)^2(2n+1)}$ and
\[
M^x_n(y) =M_n(x,y) = 1 + \frac{t_n}{t_{n-1}}xy = 1 + \frac{(2n-1)^2}{2n(2n+1)}xy
\]
i.e. $M_n = \bihomographic{(2n-1)^2}{0}{0}{2n(2n+1)}{0}{0}{0}{2n(2n+1)}$.
Let $w=x^2$ and define
\begin{equation}\label{eq:asinSeries}
a_n(w) = M^w_n(M^w_{n+1}(\cdots)) = 1 + \frac{t_n}{t_{n-1}}w+ \frac{t_{n+1}}{t_{n-1}}w^2 + \cdots
\end{equation}
so $\arcsin(x) = x a_1(w)$ and $a_n(w) = M^w_n(a_{n+1}(w))$. Since $0 \leq w \leq 1$ and $\{t_n\}$ is strictly decreasing,
\begin{equation}\label{eq:asinBound}
1 + \frac{(2n-1)^2}{2n(2n+1)}w \leq a_n(w) \leq 1 +  \frac{(2n-1)^2w}{2n(2n+1)(1-w)} \ .
\end{equation}
The lower bound is always between 1 and 2. If the upper bound is less than 2, let $[p_0,p_1\cdots p_k] = [0,p_1\cdots p_k]$ be the longest common prefix of the CL expansions of the upper and lower bounds; this common prefix will also be a prefix of $a_n(w)$. 

 The rest of $a_n(w)$ is obtained by applying the arithmetic algorithm to
\[
\hat{M_n}(w, a_{n+1}(w)) 
\]
where $\hat{M_n}$ is obtained by starting with $M_n$ and producing the terms of the longest common prefix, i.e.
\[
\hat{M_n} = \mbox{prod}_{p_k}( \cdots  \mbox{prod}_{p_0}(M_n) \cdots)
\]

If the upper bound is greater than 2, then (\ref{eq:asinBound}) does not provide enough information about $a_n(w)$, so we must use more terms of (\ref{eq:asinSeries}), enough to obtain a non-empty common prefix; in general
\[
1 + \frac{t_n}{t_{n-1}}w + \cdots + \frac{t_{n+k-1}}{t_{n-1}}w^k  < 
a_n(w) <
1 + \frac{t_n}{t_{n-1}}w + \cdots + \frac{t_{n+k-1}w^k}{t_{n-1}(1-w)}\ .
\]

\subsection{Trigonometric Functions}
The cosine
\[
\cos(x) = 1 - \frac{x^2}{2!} + \frac{x^4}{4!} - \cdots = 1 - \frac{x^2}{2} \left(1 - \frac{x^2}{12}\left( 1 - \frac{x^2}{30}\left(1 + \cdots \right. \right.  \right.  
\]
can be represented by letting $w=x^2$, defining
\[
M^w_n(y) = 1 - \frac{wy}{2n(2n-1)}
\]
and
\begin{equation}\label{eq:cosIter}
c_n(w) = M^w_n(M^w_{n+1}(\cdots)) = 1 - \frac{w}{2n(2n-1)} + \frac{w^2}{(2n+2)(2n+1)2n(2n-1)} - \cdots \ .
\end{equation}
so $c_n(w) = M^w_n(c_{n+1}(w))$ and $c_1(w) = \cos(x)$.

If $|x| \leq \frac{\pi}{2}$, then $w < \frac{5}{2}$, and (\ref{eq:cosIter}) is an alternating series bounded between $1 - \frac{5}{4n(2n-1)}$ and 1. The CL representation of $c_n(w)$ will start with $-1$; for $n \geq 2$, we then have
\[
1 < \frac{1}{c_n(w)} \leq 1 + \frac{5}{8n^2 - 4n - 5} < 2
\]
so the next term is zero, followed by $\lfloor \log_2(\frac{8n^2 - 4n - 5}{5}) \rfloor$ ones. As with $e^x$, we get subsequent terms by applying $\produce_{-1}, \produce_0$ and $\produce_1$ to $M_n$, then using the arithmetic algorithm.

If $\frac{\pi}{2} \leq |x| \leq \pi$ we can use $\cos(x) = -\cos(\pi - x)$, and for all other values $\cos(x) = \cos(x \pm 2\pi)$. We can compute $\pi$ as $6\arcsin(\frac{1}{2})$. We then define $\sin(x) = \cos(x-\pi/2)$ and $\tan(x)=\sin(x)/\cos(x)$.

\subsection{Logarithms}
To compute $\log(x)$ we use the series
\[
\log (x) =  2 \left( \left(\frac{x-1}{x+1}\right) + \frac{1}{3}\left(\frac{x-1}{x+1}\right)^3 + \frac{1}{5}\left(\frac{x-1}{x+1}\right)^5 \cdots + \right)
\]
which converges for all $x > 0$. Let $z=\frac{x-1}{x+1}$, and define
\[
g(w) = 1 + \frac{w}{3} + \frac{w^2}{5} + \frac{w^3}{7} + \cdots = 1 + \frac{w}{3}\left(1 + \frac{3w}{5}\left(1 + \frac{5w}{7}\left(1 +\cdots \right. \right. \right.\ .
\]
Then $\log(x) = 2zg(z^2)$; since $|z| \leq 1$, we only need consider how to compute $g(w)$ with $0 \leq w \leq 1$.
The series for $g(w)$ can be obtained by iterating
\[
M^w_n(y) = 1 + \frac{2n-1}{2n+1}wy 
\]
i.e. the bihomographic expression $\bihomographic{2n-1}{0}{0}{2n+1}{0}{0}{0}{2n+1}$.
Define
\[
g_n(w) = 1  + \sum_{k=1}^{\infty}\frac{2n-1}{2n+2k-1}w^k= \lim_{m\to\infty} M_{(n,w)}( M_{(n+1,w)}( \cdots M_{(m,w)}(1) \cdots ) ) \ .
\]
Then $g_1(w) = g(w)$, and $g_n(w) = M^w_n(g_{n+1}(w))$. Since $w$ is nonnegative, we have
\begin{equation}\label{eq:logbound}
1 \leq g_n(w) \leq 1 + \frac{2n-1}{2n+1}\sum_{k=1}^{\infty} w^k = 1+\frac{2n-1}{2n+1}\frac{w}{1-w}\ .
\end{equation}
If $1 \leq x \leq 2$, then $\frac{w}{1-w} \leq \frac{1}{8}$, so $g_n(w) \leq 2$ and the first term of its CL is zero. Then
\[
\frac{1}{1-g_n(w)} \geq \frac{8(2n+1)}{2n-1}
\]
and the CL continues with $\lfloor \log_2(\frac{8(2n+1)}{2n-1}) \rfloor$ ones et cetera.

For $x<1$ we can use $\log(x) = -\log(1/x)$, and for $x>2$ we can use $\log(x) = \log(2) + \log(x/2)$.

\section{Termination}\label{sec:termination}
As noted previously, the arithmetic algorithm is not guaranteed to produce an output term after reading a finite amount of input. Arithmetic on continued fractions has precisely the same problem. In \cite{collins2025transcendentalfunctionscontinuedfractions}, we addressed the problem by extending the internal representation of the continued fraction data type to include upper and lower bounds on $M(x,y)$; we include these bounds as part of the output, and nested arithmetic can use such bounds on $x$ and $y$ to determine the next bound on $M(x,y)$. When we obtain $n - \varepsilon < M(x,y) < n + \varepsilon$, and $\varepsilon$ is less than a specified error threshold, we consider the next term of the CF to be $n$. But we retain the bounds, in case subsequent computations require greater precision.

Such an approach is less appealing with continued logarithms, since part of the motivation for studying CLs is their amenability to very compact hardware implementation; manipulating arbitrary rational bounds adds considerable complexity to the prod/con manipulations defined in section \ref{sec:arithAlg}. One option might be to have a parallel error channel which detects and reports stalled computations.

\bibliographystyle{plain}
\bibliography{clogs}

\end{document}